\documentclass{cpaper}

\usepackage{amsfonts}
\usepackage[dvips]{graphics,color}
\newcommand{\Z}{\mathbb{Z}}

\input amssym.def
\newsymbol\square 1003
\newsymbol\lozenge 1006
\newsymbol\surjarrow 1310
\newsymbol\injarrow 131A
\newsymbol\ltimes 226E
\newsymbol\rtimes 226F
\newsymbol\void 203F
\newsymbol\leftrightsquigarrow 1321

\def\proclaim #1. #2\par #3\par {\medbreak
\noindent{\bf#1.\enspace}{\sl#2}\par\medbreak
\noindent{\bf Proof.} #3 \par}
\def\proclaimb #1. #2\par {\medbreak
\noindent{\bf#1.\enspace}{\sl#2}\par\medbreak}
 0
 2
 0
\newcommand{\arabiceqn}{\setcounter{equation}{0}%
\renewcommand{\theequation}{\mbox{\arabic{section}.\arabic{equation}}}}%
\newtheorem{thm}[equation]{Theorem}
\newtheorem{lemma}[equation]{Lemma}
\newtheorem{prop}[equation]{Proposition}

\newtheorem{claim}[equation]{Claim}

\begin{document}
\title{Contractibility of fixed point sets of auter space}
\author{Craig A. Jensen}
\date{Spring 2000}

\smallskip

\address{Department of Mathematics, The Ohio State University\\
Columbus, OH 43210, USA}
\email{jensen@math.ohio-state.edu}

\begin{abstract}
We show that for every finite subgroup $G$ of $Aut(F_n)$, the fixed point 
subcomplex $X_n^G$ is contractible, where $F_n$ is the free group 
on $n$ letters and $X_n$ is the spine of 
``auter space'' constructed by Hatcher and Vogtmann in \cite{[H-V]}.
In more categorical language, $X_n =\underline{E}Aut(F_n)$.
This is useful because it allows one to compute	
(see, for example, \cite{[JD],[J1]})
the cohomology of normalizers or centralizers of finite subgroups
of $Aut(F_n)$ based on their actions on fixed point subcomplexes.
The techniques used to prove it are largely those of Krstic and Vogtmann
in \cite{[K-V]}, who in turn used techniques similar to Culler and 
Vogtmann in \cite{[C-V]}
\end{abstract}


\maketitle

\section{Introduction}
\arabiceqn

Let $F_n$ denote the free group on $n$ letters and let
$Aut(F_n)$ and $Out(F_n)$ denote the automorphism group
and outer automorphism group, respectively, of $F_n$.
In \cite{[C-V]} Culler and Vogtmann defined a space on which
$Out(F_n)$ acts nicely called ``outer space''.  By
studying the action of $Out(F_n)$ on this space, various
people have been able to calculate the cohomology
of $Out(F_n)$ in specific cases.  More recently, Hatcher in
\cite{[H]} and Hatcher and Vogtmann in \cite{[H-V]} have defined
a space on which $Aut(F_n)$ acts nicely called ``auter space''
and have used this to calculate the cohomology of
$Aut(F_n)$ in specific cases.

We review some basic properties and definitions of 
auter space.  Most of these can be found in \cite{[C-V]},
\cite{[H-V]}, \cite{[S-V1]}, or \cite{[S-V2]}.
Let $(R_n,v_0)$ be the
$n$-leafed rose,
a wedge of $n$ circles. We say a pointed graph $(G,x_0)$
is \emph{admissible} if it has no free edges, all vertices except the
basepoint
have valence at least three, and there is a basepoint-preserving continuous
map $\phi\colon R_n \to G$ which induces an isomorphism on $\pi_1$. The
triple $(\phi,G,x_0)$ is called a \emph{marked graph}.  Two marked graphs
$(\phi_i,G_i,x_i) \hbox{ for } i=0,1$  are \emph{equivalent}  if there is
a homeomorphism $\alpha \colon (G_0,x_0) \to (G_1,x_1)$ such that
$ (\alpha\circ\phi_0)_\# = (\phi_1)_\# : \pi_1(R_n,v_0) \to \pi_1(G_1,x_1)$.
Define a partial order on the set of all equivalence classes of
marked graphs by setting $(\phi_0,G_0,x_0) \leq
(\phi_1,G_1,x_1)$ if $G_1$ contains a \emph{forest}
(a disjoint union of trees in $G_1$ which contains all of the vertices
of $G_1$) such that collapsing each tree in
the forest to a point yields $G_0$, where the collapse is compatible with
the maps $\phi_0$ and $\phi_1$.

From \cite{[H]} and \cite{[H-V]} we have that $Aut(F_n)$ acts
with finite stabilizers on
a contractible space $X_n$.
The space $X_n$ is the geometric realization of the poset of
marked graphs that we defined above.
Let $Q_n$ be the quotient of $X_n$ by
$Aut(F_n)$.
Note that the CW-complex $Q_n$ is not necessarily a simplicial
complex.
Since $Aut(F_n)$ has a torsion free subgroup of finite index
\cite{[H]} and it acts on the contractible,
finite dimensional space
$X_n$ with finite stabilizers and finite quotient,
$Aut(F_n)$ has finite vcd. From \cite{[Z]} (cf. \cite{[C]}), any finite 
subgroup $G$ of $Aut(F_n)$ fixes a point of $X_n$. 
Our goal is to show

\begin{thm} \label{tr13} Auter space is an $\underline{E}Aut(F_n)$-space.
That is, for any finite subgroup $G$ of $Aut(F_n)$,
the fixed point subcomplex $X_n^G$ is contractible.
\end{thm}

This paper is based in part on a dissertation written while the author 
was a student of Karen Vogtmann at Cornell, and the author would like
to thank Prof. Vogtmann for her help and advice.

\section{Norms and Absolute Values} \label{c11}
\arabiceqn
We strongly recommend that the reader study \cite{[K-V]} by 
Krstic and Vogtmann, where they prove the analog of 
Theorem \ref{tr13} for $Out(F_n)$ and outer space.
This paper is essentially a modification of their results on
fixed point spaces of outer space to fixed point spaces
of auter space, and we will often omit details which are
similar to work already done in \cite{[K-V]}.
White \cite{[W]} also proved the
result for fixed point subcomplexes of outer space, 
but we do not know to 
what extent his work can be applied to auter space.

In particular, Krstic and Vogtmann define a complex $L_G$ of
``essential marked $G$-graphs'' that the fixed point set
$X_n^G$ in outer space deformation retracts to.  Then they
order the reduced marked $G$-graphs in $L_G$ using
a norm $\|\cdot\|_{out}$.  Using this norm
to determine which reduced marked $G$-graphs
should be considered next, Krstic and
Vogtmann performed a transfinite induction argument to
show that $L_G$ is contractible, by
building $L_G$ up as the union of stars
of reduced marked $G$-graphs.

We will follow a similar approach, and define norms
$$\|\cdot\|_{aut}
\hbox{ and }
\|\cdot\|_{tot}=\|\cdot\|_{out} \times \|\cdot\|_{aut}$$
to order the reduced marked essential $G$-graphs
in auter space.  For technical reasons,
$\|\cdot\|_{tot}$ will be the appropriate
norm to use when performing the transfinite
induction argument to show the contractibility
of the corresponding $L_G$ in auter space.

The norm $\|\cdot\|_{out}$ was defined by Krstic and
Vogtmann as follows.  Order the set $\mathcal{W}$ of
conjugacy classes of elements of $F_n$ as
$\mathcal{W} = \{w_1, w_2, \ldots \}.$
Totally order $\Z^\mathcal{W}$ by the lexicographic order.
 Let
$\sigma = [s, \Gamma]$ be a marked graph and
define $\|\sigma\|_{out} \in \Z^\mathcal{W}$
by letting $(\|\sigma\|_{out})_i$ be the sum
over all $x \in G$ of the
lengths in $\Gamma$ of the reduced loops (given by the
marking $s$) corresponding to $xw_i$.  Equivalently,
they define an absolute value
$|\cdot|_{out} \in \Z^\mathcal{W}$ on the edges
of $\Gamma$ and set
$$\|\sigma\|_{out} = \frac{1}{2}\sum_{e \in E(\Gamma)} |e|_{out}.$$
The $i$th coordinate of $|e|_{out}$ is simply the
sum for all $x \in G$ of the contributions
of $e$ or $\bar e$ to the loop $xw_i$ in $\Gamma$.
In other words, it is the sum over all $x \in G$ of
the number of times $e$ or $\bar e$ appears in the
cyclically reduced edge path representing $xw_i$.
For $A,B \subseteq E(\Gamma)$ define
$(A.B)_{out} \in \Z^\mathcal{W}$
to be the function whose $i$th coordinate is
the sum over all $x \in G$ of the number of
times $a\bar b$ or $b\bar a$ appears in the
reduced loop in $\Gamma$ corresponding to $xw_i$.
Finally, for $C \subseteq E(\Gamma)$, define
$|C|_{out}$ inductively by the formula
$$\matrix{\hfill |A\coprod B|_{out} = |A|_{out} + |B|_{out} - 2(A.B)_{out} \hfill \cr}$$
for disjoint subsets $A$ and $B$ of $E(\Gamma)$.
Note that with the above definition,
$|A|_{out}=(A.E(\Gamma)-A)_{out}=|E(\Gamma)-A|_{out}$.

The corresponding quantities for $Aut(F_n)$ are defined
in much the same way, the basic difference being that
we think of the lengths of reduced paths rather than reduced
loops.  Order $F_n$ as $F_n = \{ \alpha_1, \alpha_2, \ldots \}$,
and give $\Z^{F_n}$ the lexicographic order.
For a finite subgroup $G$ of $Aut(F_n)$, consider
a pointed marked $G$-graph $\sigma = [s,\Gamma]$.
Define the norm $\|\sigma\|_{aut} \in \Z^{F_n}$
to be $|G| \cdot L$, where $L : F_n \to \Z$ is
the Lyndon length function of the marked graph.
In other words, the pointed marked graph $\sigma$ corresponds
to an action of $F_n$ on a rooted $\Z$-tree $T$.
Define $$L(\alpha_i) = \{\hbox{the distance } \alpha_i  \hbox{ moves
the root of } T \}.$$
Equivalently, the $i$th coordinate of $\|\sigma\|_{aut}$
is the sum over all $x \in G$ of the lengths in
$\Gamma$ of the reduced (but {\em not} cyclically reduced)
paths corresponding to $x\alpha_i \in \pi_1(\Gamma,*)$.

As before in the case of $Out(F_n)$, we can
define an absolute value
$|\cdot|_{aut} \in \Z^{F_n}$ on the edges
of $\Gamma$ and set
$$\|\sigma\|_{aut} = \frac{1}{2}\sum_{e \in E(\Gamma)} |e|_{aut}.$$
The $i$th coordinate of $|e|_{aut}$ is simply the
sum of for all $x \in G$ of the contributions
of $e$ or $\bar e$ to the reduced
(but not cyclically reduced) path $x\alpha_i$ in $\pi_1(\Gamma,*)$.
Hence it is the sum over all $x \in G$ of
the number of times $e$ or $\bar e$ appears in the
reduced edge path representing $x\alpha_i$.
For $A,B \subseteq E(\Gamma)$ define
$(A.B)_{aut} \in \Z^\mathcal{W}$
to be the function whose $i$th coordinate is
the sum over all $x \in G$ of the number of
times $a\bar b$ or $b\bar a$ appears in the
reduced path in $\Gamma$ corresponding to $x\alpha_i$.
Finally, for $C \subseteq E(\Gamma)$, define
$|C|_{aut}$ inductively by the formula
$$\matrix{\hfill |A\coprod B|_{aut} = |A|_{aut} + |B|_{aut} - 2(A.B)_{aut} \hfill \cr}$$
for disjoint subsets $A$ and $B$ of $E(\Gamma)$.
In contrast to the case with $Out(F_n)$ the formula
$|A|_{aut}=(A.E(\Gamma)-A)_{aut}$
certainly does not hold any longer.

Our final norm $\|\cdot\|_{tot}$ is just the
product of the previous two.  That is, let
$\sigma = [s,\Gamma]$ be a pointed marked
$G$-graph for a finite subset $G$ of $Aut(F_n)$
and totally order $\Z^\mathcal{W} \times \Z^{F_n}$
by the lexicographic order.  Define
$\|\sigma\|_{tot} \in \Z^\mathcal{W} \times \Z^{F_n}$
as $\|\sigma\|_{tot} = \|\sigma\|_{out} \times \|\sigma\|_{aut}$,
where to calculate $\|\sigma\|_{out}$ we just
forget that $\Gamma$ has a basepoint.
The functions $|e|_{tot}$, $(A.B)_{tot}$, and
$|A|_{tot}$ are defined similarly.

For a vertex $v$, let $E_v$ be the set of oriented
edges ending at $v$.  
We call certain subsets $\alpha \subseteq E_v$ {\em ideal edges}
and think of them as corresponding to new edges created when we blow up
the original graph at the vertex $v$ by pulling away the edges in $\alpha$.
Formally,
the notion of ideal edges is defined as in \cite{[K-V]},
with the exception that if the ideal edge
$\alpha  \subseteq E_*$ then condition $(i)$ of their definition
should be changed to:
$$(i) \hbox{ } card(\alpha) \geq 2 \hbox{ and } card(E_* - \alpha) \geq 1.$$
That is, ideal edges at the basepoint can contain all except
one of the edges of $E_*$.  The definition of blowing up an
ideal edge is taken exactly as defined in \cite{[K-V]}.  Hence
if we are blowing up an ideal edge $\alpha \subseteq E_*$ then
we are pulling the edges of $\alpha$ away from the basepoint
along a new edge $e(\alpha)$ we just constructed.  If
$card(E_* - \alpha) = 1$, this will result in a graph whose
basepoint has valence $2$.

Let $\alpha$ be an ideal edge of $\sigma = [s, \Gamma]$ and
$\sigma^{G\alpha} = [s^{G\alpha}, \Gamma^{G\alpha}]$
be the result of blowing up the ideal edge $\alpha$.
Then it is easy to show that $|\alpha|_{aut}$ in
$\Gamma$ is equal to $|e(\alpha)|_{aut}$ in $\Gamma^{G\alpha}$
(which was the whole point of defining $|\cdot|_{aut}$
on subsets of edges.)  Hence $|\alpha|_{tot}= |e(\alpha)|_{tot}$
also, as Krstic and Vogtmann show the corresponding formula
for $|\cdot|_{out}$.  From this, the analogs of Proposition
6.4  about Whitehead moves in \cite{[K-V]} are true
for the norms $\|\cdot\|_{aut}$ and
$\|\cdot\|_{tot}$.  That is, for an ideal edge $\alpha$ define
$D(\alpha)$ by
$$D(\alpha) = \{ a \in \alpha : stab(a)=stab(\alpha) \hbox{ and }
\bar a \not \in \bigcup G \alpha \}.$$
Then the Whitehead move $(G\alpha,Ga)$ is the result
of first blowing up $\alpha$ in $\Gamma$ to get $\sigma^{G\alpha}$ and
then collapsing $Ga$ in $\Gamma^{G\alpha}$ to get
$\sigma'$.  Proposition 6.4 of \cite{[K-V]} states that
$$\|\sigma'\|_{out} = \|\sigma\|_{out} + [G:stab(\alpha)]
(|\alpha|_{out} - |a|_{out}).$$
As mentioned before, this remains true if out-norms and
absolute values are replaced by aut- or tot-norms and
absolute values.

The value $[G:stab(\alpha)] (|a|_{out} - |\alpha|_{out})$
is called the {\em out-reductivity} of $(\alpha,a)$
and is denoted $red_{out}(\alpha,a)$.  Similar
notions of {\em aut-reductivity} and {\em tot-reductivity}
are defined as well.  A Whitehead move reduces the norm
iff the corresponding reductivity is greater than zero,
in which case the Whitehead move is called
{\em reductive}.
The {\em $x$-reductivity} of
an ideal edge $\alpha$ is the maximum over all
elements $a \in D(\alpha)$ of $red_x(\alpha,a)$,
where $x$ is out, aut, or tot.  It thus makes
sense to talk of an ideal edge $\alpha$ as being
{\em out-reductive}, etc. The norm $\|\cdot\|_{tot}$
will be useful to us because:

\begin{prop} \label{t17} Let $\alpha \subseteq E_v$
be a tot-reductive ideal edge of a reduced
marked $G$-graph $\rho$.  Suppose $\alpha$ is {\em invertible}
(that is, $E_v - \alpha \not \subseteq G\alpha$
and $E_v - \alpha$ is an ideal edge.)  Then
$\alpha^{-1} = E_v - \alpha$ is tot-reductive.
\end{prop}

\begin{proof} Assume $v=*$, else the proof is trivial.  Say $(\alpha,a)$ is the
reductive ideal edge.  Since $stab(*)=G$ and $\alpha$ is
invertible, the analog of Lemma 5.1 of \cite{[K-V]}
gives us that $stab(\alpha)=stab(a)=G$.  Say $\rho=[s,\Gamma].$
As before, let $\rho^{G\alpha}=[s^{G\alpha},\Gamma^{G\alpha}]$
be the result of blowing up the ideal edge $(\alpha,a)$.
Then let $\rho'=[s',\Gamma']$ be the result of collapsing
$a$ in $\Gamma^{G\alpha}$.  We know that
$\|\rho'\|_{tot} < \|\rho\|_{tot}$ as $(\alpha,a)$
is tot-reductive.

Assuming the claim below, it will be easy to complete the
proof as follows:  Since $\|\rho'\|_{out} \not = \|\rho\|_{out}$
and $\|\rho'\|_{tot} < \|\rho\|_{tot}$,
we must have $\|\rho'\|_{out} < \|\rho\|_{out}$.
Let $\rho''$ be the result of doing the Whitehead move
$(\alpha^{-1},a^{-1})$ to $\rho$.  Because
$red_{out}(\alpha,a)=red_{out}(\alpha^{-1},a^{-1})$
(see the comments in \cite{[K-V]} following the proof
of \S 6.4), it follows that
$\|\rho''\|_{out} = \|\rho'\|_{out}$.
So $\|\rho''\|_{out} < \|\rho\|_{out}$ and
hence $\|\rho''\|_{tot} < \|\rho\|_{tot}$.
Thus $\alpha^{-1}$ is reductive. \end{proof}

\begin{claim} \label{t18} $\|\rho'\|_{out} \not = \|\rho\|_{out}$.
\end{claim}

\begin{proof} Since $stab(a)=G$ and $\rho = [s,\Gamma]$ is reduced, the edge $a$
must both begin and end at $*$.
Enumerate the edges of $\alpha - \{a\}$ and $\alpha^{-1} - \{a^{-1}\}$
as $b_0, \ldots, b_r$ and $c_0, \ldots, c_s$, respectively,
where $r,s \geq 0$. 
We have three cases, which are not disjoint but
are exhaustive.
\begin{enumerate}
\item Some $b_i$ is a loop at $*$ and $b_i^{-1} \not \in \alpha$.
Let $w_k \in \mathcal{W}$ be an element that maps to the
loop $b_i$.  Then $(|a|_{out})_k = 0$ and
$(|\alpha|_{out})_k \geq 1$ since the loop $b_i$ is
sent to $b_ie(\alpha)$.

\item Some $b_i$ starts at another vertex $v_1 \not = *$.
Since $G$ acts nontrivially on $b_i$ and because $b_i$
must be elliptic (as it is clearly not bent hyperbolic),
there must be another $b_j \not = b_i$
also going from $*$ to $v_1$.
(For the definitions
of elliptic and bent hyperbolic see \S 4A in the paper
by Krstic and Vogtmann.)  There are two subcases:
\begin{itemize}
\item There is an edge $c_l$ in $\alpha^{-1} - \{a^{-1}\}$ that
begins and ends at $*$.  We can assume
$c_l^{-1} \in \alpha^{-1} - \{a^{-1}\}$ also,
else we are in case 1.  Choose a $w_k \in \mathcal{W}$ that
maps to the loop $b_ib_j^{-1}c_l$.
Now $(|a|_{out})_k = 0$ and
$(|\alpha|_{out})_k \geq 1$ since $w_k$ is
sent to $e(\alpha)^{-1}b_ib_j^{-1}e(\alpha)c_l$.
\item There is an edge $c_l$ in $\alpha^{-1} - \{a^{-1}\}$ that
begins at $v_2 \not = *$ and ends at $*$.
Because $G$ acts nontrivially on $c_l$ and $c_l$ is
elliptic, there is another
edge $c_m \not = c_l$ also going from $v_2$ to $*$.
Choose a $w_k \in \mathcal{W}$ that
maps to the loop $b_ib_j^{-1}c_lc_m^{-1}$.
Then $(|a|_{out})_k = 0$ but
$(|\alpha|_{out})_k \geq 1$ as $b_ib_j^{-1}c_lc_m^{-1}$ is
sent to $e(\alpha)^{-1}b_ib_j^{-1}e(\alpha)c_lc_m^{-1}$.
\end{itemize}

\item Some $b_i$ is a loop at $*$ and $b_i^{-1} \not \in \alpha$
also.  As in case $2.$ above, there are two subcases.
\begin{itemize}
\item Same as in case $2.$ above.  Choose a $w_k \in \mathcal{W}$
that maps to $b_ic_l$.
Then $(|a|_{out})_k = 0$ but
$(|\alpha|_{out})_k \geq 1$ as $b_ic_l$ is
sent to $e(\alpha)^{-1}b_ie(\alpha)c_l$.
\item Same as in case $2.$ above.  Choose a $w_k \in \mathcal{W}$
that maps to $b_ic_lc_m^{-1}$.
Hence $(|a|_{out})_k = 0$ and yet
$(|\alpha|_{out})_k \geq 1$ because $b_ic_lc_m^{-1}$ is
mapped to $e(\alpha)^{-1}b_ie(\alpha)c_lc_m^{-1}$.
\end{itemize}
\end{enumerate}

In each case we have $|a|_{out} \not = |\alpha|_{out}$;
therefore, $red_{out}(\alpha,a) \not = 0$
and $\|\rho'\|_{out} = \|\rho\|_{out}$. \end{proof}

Because of Proposition \ref{t17}, tot-reductivity
will be the most useful of the three types of reductivity
(out, aut, and tot) for us.
From now on when we say that $\rho$ is
{\em reductive}, this is just shorthand for saying $\rho$ is
tot-reductive.

Proposition $6.1$ of \cite{[K-V]}, states that
$$\matrix{\hfill |A\coprod B|_{out} = |A|_{out} + |B|_{out} - 2(A.B)_{out} \hfill \cr}$$
for disjoint subsets $A$ and $B$ of $E(\Gamma)$.
This also holds for $aut$-norms because it
is our definition of the absolute values $|\cdot|_{aut}$
for sets of edges and can be inductively shown
to be well-defined.
It is important that this
property holds for $aut$-norms because
it is used by
many of the later propositions in Krstic and
Vogtmann (e.g., Proposition $6.2$ of \cite{[K-V]}
which will correspond to our Proposition \ref{t19}.)

Proposition $6.2$ of \cite{[K-V]} states that:

\begin{prop}[Krstic-Vogtmann] \label{tr38}
Let $K$ be a subgroup of $G$, let $A$ be a $K$-invariant
subset of $E(\Gamma)$, and let $e$ be an edge of $\Gamma$
with $stab(e)$ contained in $K$.  Then
$$((Ke).A)_{out} = [K:stab(e)](e.A)_{out}.$$
\end{prop}

We now show Proposition $6.2$
of \cite{[K-V]} also holds
for the aut-norm,
which will be useful in some
combinatorial lemmas later in this section.
Once we show that the analog of Proposition \ref{tr38}
is true for the $aut$-norm,
it
will be true for both the out- and aut-norms
on a component-by-component basis.  In other words,
the equality stated in the proposition
is true for each component of
$\Z^\mathcal{W}$ or $\Z^{F_n}$ and does not use the
total (lexicographic) order on those sets.  Hence
it is automatically true for the tot-norm, as the
tot-norm is just the product of the out-norm
and the aut-norm.  We will be able to use the
same approach (that of just showing something to be
true for the aut-norm)
in some lemmas later on in this section.

\begin{prop} \label{t19}
Let $K$ be a subgroup of $G$,
$A$ be a $K$-invariant subset of $E(\Gamma)$,
and $e$ be an edge of $\Gamma$ with $stab(e)$
contained in $K$.  Then
$$((Ke).A)_{aut}=[K:stab(e)](e.A)_{aut}.$$
\end{prop}

\begin{proof} To simplify the notation in the proof below,
we write (just for this proof)
$\|\cdot\|$ for $\|\cdot\|_{aut}$,
$|\cdot|$ for $|\cdot|_{aut}$,
reductive for aut-reductive, etc.
 
Examine $((Ke).A)_i$.  It is the number of
times one of the strings $(ke)a^{-1}$ or
$a(ke)^{-1}$ appears in one of the
$x\alpha_i$, for all $k \in K$, $a \in A$, and
$x \in G$.

Now $stab(e) \subseteq K$ and we can write
$$\matrix{\hfill K = stab(e) \coprod k_2stab(e) \coprod \ldots
\coprod k_{[K:stab(e)]}stab(e) \hfill \cr}$$
using coset representatives $k_i$.
Further note that the number of times one
of the strings $ea^{-1}$ or $ae^{-1}$
appears in one of the strings $x\alpha_i$
for $a \in A$, $x \in G$ is exactly the same
as the number of times one of the strings
$k_iea^{-1}$ or $a(k_ie)^{-1}$ appears in
one of the $x\alpha_i$
for $a \in A$, $x \in G$.  This is because
each $k_i$ is in $G$ and $A$ is
$K$-invariant so if $ea^{-1}$ is in $x\alpha_i$
then $(k_ie)(k_ia)^{-1}$ is in
$(k_ix)\alpha_i$.
So $((Ke).A)_i = [K:stab(e)] (e.A)_i$. \end{proof}

\begin{prop} \label{t20} The set of pointed marked
$G$-graphs is well-ordered by the tot-norm.
\end{prop}

\begin{proof} Let $\mathcal{A}$ be a nonempty collection of pointed
marked $G$-graphs.  We must find a least element of $\mathcal{A}$.
Let $[\mathcal{A}]$ be the set of equivalence classes of
marked $G$-graphs in $\mathcal{A}$ obtained by forgetting
the basepoint $*$.  From
Proposition $6.3$ of \cite{[K-V]} the out-norm
well orders marked $G$-graphs, and  $[\mathcal{A}]$ has a
least element $U \subseteq \mathcal{A}$.

Say $\sigma = [s,\Gamma]$ is the marked
$G$-graph representing this $U$.  The marked graph $\sigma$
corresponds to an action of $F_n$ on the tree
$\tilde \Gamma = \Lambda$.
From \cite{[C-V]} $\sigma$ corresponds to a free, {\em minimal} (there
are no invariant proper subtrees), and not abelian
(an action is {\em abelian} iff every element of the
commutator $[F_n,F_n]$ has length $0$) action without
inversions on the tree $\tilde\Gamma = \Lambda$. 

The action has an associated non-abelian (see Alperin and Bass
in \cite{[A-B]}) length function $l$ on $F_n$.
By Theorem $7.4$ of \cite{[A-B]}, there exist hyperbolic
elements $\alpha_n, \alpha_m \in F_n$, $n < m$,  
such that the characteristic subtrees $A_{\alpha_n}$ and 
$A_{\alpha_m}$ are linear and disjoint.

Recall that we wish to find the least element of $U$ in the 
tot-norm.  Following the
proof of Proposition $6.3$ in \cite{[K-V]}, we set
$U_0 = U$ and define $U_i$ inductively for $i \geq 1$.
Let $\gamma_i = min\{(\|\delta\|_{aut})_i : \delta \in U_{i-1}\}.$
Next define $U_i$ to be the subset of $U_{i-1}$ consisting of
$\delta$ with $(\|\delta\|_{aut})_i = \gamma_i$.
To finish our proof, it suffices to show that $U_m$ has only
finitely many elements.

Each element of $U$ corresponds to an action of $F_n$
on a pointed tree.  In each case, if we forget the
basepoint then the tree is homeomorphic to $\Lambda$.  
The map from $U$ to Lyndon length functions on $F_n$, given by
seeing how far the basepoint is moved under the corresponding
action, is injective (see \cite{[H-V]}, \cite{[A-B]}.)  
Note that in each case, the action of $F_n$
on the underlying non-pointed tree $\Lambda$ is the
same.  We are only varying where we place the basepoint on $\Lambda$
and seeing how far elements of $F_n$ move this
basepoint.

The elements of $U_1$ are those where the basepoint is located
closest to the linear subtree $A_{\alpha_1} \subset \Lambda$, and
$U_1$ could be infinite. Let $B$ be the bridge joining $A_{\alpha_n}$
and $A_{\alpha_m}$. 
To show that $U_m$ is finite,
 it suffices to show that there are only
finitely many points at fixed distances $d_1$ and
$d_2$ from $A_{\alpha_n}$ and $A_{\alpha_m}$, respectively.
If $d(x,A_{\alpha_n})=d_1$ and $d(x,A_{\alpha_m})=d_2$, then
choose paths $p_1$ and $p_2$ of lengths $d_1$ and $d_2$
from $x$ to $q_1 \in A_{\alpha_n}$ and
$q_2 \in A_{\alpha_m}$, respectively.  The union of these
two paths $p_1$ and $p_2$ contains the bridge $B$.
Consequently, $d(x,B) \leq d_1 + d_2$.  Since the tree
is locally finite and $B$ is finite, $x$ is one of a finite number
of vertices. \end{proof}

A few definitions are in order at this point.  Basically,
we are trying to find the appropriate parallels of
definitions in \cite{[K-V]}.  Fix a reduced marked
$G$-graph $\rho = [s,\Gamma]$.  Let $(\mu,m)$ be
a maximally reductive ideal pair of $\rho$.  That is,
$\mu$ is the maximally reductive ideal edge in $\rho$
and $m \in D(\mu)$ is an edge in $\mu$ which
allows the Whitehead move $(\mu,m)$ to realize this
maximum.

Let $\alpha \subset E_u$ and $\beta \subset E_v$ be ideal
edges of $\rho$.  Then the ideal edge orbits $G\alpha$
and $G\beta$ are {\em compatible} if one of the
following holds:
\begin{enumerate}
\item $G\alpha \subseteq G\beta.$
\item $G\beta \subseteq G\alpha.$
\item $G\alpha \cap G\beta = \void$ and $\alpha \not = \beta^{-1}.$
\item $G\alpha \cap G\beta = \void$ and $u=v=*.$
\end{enumerate}
The ideal edge orbits $G\alpha$
and $G\beta$ are {\em pre-compatible} if one of the
following holds:
\begin{enumerate}
\item They are compatible.
\item $\alpha$ is invertible and $\alpha^{-1} \subseteq \beta$.
\item $\beta$ is invertible and $\beta^{-1} \subseteq \alpha$.
\end{enumerate}
Note that $2.$ and $3.$ above would be equivalent if we did not
need to consider ideal edges of the form $\gamma = E_* - \{c^{-1}\}$
which have $stab(\gamma)=G$ but are not invertible.

An {\em oriented ideal forest} is a collection of
pairwise compatible ideal edge orbits.  These can be blown up
to obtain marked graphs in the star in $L_G$ of $\rho$.  The
correspondence is not unique, however, as two different
oriented ideal forests can be blown up to yield the same
marked graph.  This problem is solved by defining ideal
forests.  There is a poset isomorphism between the poset
of ideal forests and the star of $\rho$ in $L_G$.

An {\em ideal forest} is a collection $\Phi = \Phi_1 \coprod \Phi_2$
where $\Phi_1$ are the edges at $*$ and $\Phi_2$ are the edges not at $*$,
such that
\begin{enumerate}
\item The elements of $\Phi_2$ are pairwise pre-compatible
and $\Phi_2$ contains the inverse of each of its invertible
edge orbits; and
\item The elements of $\Phi_1$ are pairwise compatible.
\end{enumerate}

With respect to a particular reduced marked $G$-graph $\rho$
and maximally reductive ideal edge $(\mu,m)$, the following definitions
will be used frequently in the next section (which contains the core
proof of the contractibility of $L_G$.)
\begin{itemize}
\item $\mathcal{R} = \{\hbox{reductive ideal edges}\}$.
\item If $\mathcal{C}$ is a set of ideal edges, then let
$\mathcal{C}^\pm$ denote the set obtained by adjoining
to $\mathcal{C}$ the inverses of its invertible
elements that are not at the basepoint.
\item Let $S(\mathcal{C})$ be the subcomplex of
the star $st(\rho)$ spanned by ideal forests of
$\rho$, all of whose edges are in $\mathcal{C}$.
Note:  The empty forest should \underline{not}
be taken to be in $S(\mathcal{C})$.
\item $\mathcal{C}_0 = \{\alpha \in \mathcal{R} :
\alpha \hbox{ is compatible with } \mu\}$.
\item $\mathcal{C}_0^{'} = \mathcal{C}_0 \cup \{\alpha \in \mathcal{R} :
\hbox{ if } \alpha \subset E_v \hbox{ then }
stab(\alpha)=stab(v)\}$ (cf. Lemma $5.1$ of \cite{[K-V]}.)
\item $\mathcal{C}_1 = \mathcal{C}_0^{'} \cup
\{\alpha \in \mathcal{R} : m \in G\alpha \hbox{ and }
N(G\alpha,G\mu)=1\}$.
\end{itemize}
The definition of the crossing number $N(G\alpha,G\mu)$
comes from \S 7 of \cite{[K-V]} where it and other
combinatorial notions are defined.  For the
reader's convenience, we briefly state their
definitions again here.  Say
$\alpha$ and $\beta$ are two ideal edges at some
vertex $v$, with stabilizers $P$ and $Q$, respectively,
of indices $p$ and $q$ in $G$.
Choose double coset representatives
$x_1, \ldots, x_k$ of $P\backslash G/Q$. The intersection
$\delta = \alpha \cap G\beta$ breaks up as
a disjoint union
$$\matrix{\hfill \delta = \gamma_1 \coprod \ldots \coprod \gamma_k \hfill \cr}$$
with each $\gamma_i = \alpha \cap Px_i\beta$.
The $\gamma_i$ are called the {\em intersection
components} of $\alpha$ with $\beta$
and the number $N(G\alpha,G\beta)$ of nonempty
intersection components is
called the {\em crossing number}.
If $N(G\alpha,G\beta)=1$ then $G\alpha$ and $G\beta$
are said to {\em cross simply}.

The following two lemmas are stated for the out-norm
by Krstic and Vogtmann.  We will show them for
the aut-norm.
The proofs will be routine,
although they are not the same as the proofs given in
\cite{[K-V]}. This is because their proofs use the
fact that $|A|_{out} = (A.E(\Gamma)-A)_{out}$, which
is no longer true with the new norms. 
As with Proposition \ref{t19}, the
lemmas are true for both the out- and aut- norms
on a component-by-component basis.  That is,
the inequalities stated in the lemmas
are true for each component of
$\Z^\mathcal{W}$ or $\Z^{F_n}$ and do not use the
total (lexicographic) order on those sets.  Hence
it suffices to show them for the aut-norm, as the
tot-norm is the product of the out-norm
and the aut-norm.

\begin{lemma} \label{t21} 
Suppose $G\alpha$ and $G\beta$ cross simply,
with $P \leq Q$, then
$$p|\alpha \cap \beta|_{aut} + q|\beta \cup Q\alpha|_{aut}
\leq p|\alpha|_{aut} + q|\beta|_{aut}.$$
\end{lemma}

\begin{proof} To simplify the notation in the proof below,
we write (just for this proof)
$\|\cdot\|$ for $\|\cdot\|_{aut}$,
$|\cdot|$ for $|\cdot|_{aut}$,
reductive for aut-reductive, etc.

Let $[Q:P]=n$. Then $p=nq$.  Dividing by $q$,
we see that we want to show that
$$n|\alpha \cap \beta| + |\beta \cup Q\alpha|
\leq n|\alpha| + |\beta|.$$

Let $q_1, \ldots, q_n$ be a set of coset
representatives for $P$ in $Q$.
Let $\delta = \alpha \cap \beta$,
$A=\alpha - \delta$, and
$B=\beta - Q\delta$.
Since 
$$\matrix{
\hfill n|\delta| + |B \coprod Q\alpha| 
= n|\delta| + |B| + |Q\delta| + |QA| - 2 Q\delta.QA - 2 B.Q\alpha. \hfill \cr}$$
and
$$\matrix{\hfill n|\delta \coprod A| + |B \coprod Q\delta| =
n|\delta| + n|A| - 2n \delta.A + |B| + |Q\delta| - 2 B.Q\delta, \hfill \cr}$$
we have reduced the problem to showing that
$$|QA| - 2 Q\delta.QA - 2 B.Q\alpha \leq
n|A| - 2n \delta.A - 2 B.Q\delta.$$

Note that $-2 B.Q\alpha \leq -2 B.Q\delta$
as $Q\delta \subseteq Q\alpha$.
Also note that by decomposing $QA$ into a disjoint
union of $q_iA$'s, we have $Q\delta.QA \geq n\delta.A$.
Similarly, we could use induction to show that
$|QA| \leq n|A|$.  \end{proof}

\begin{lemma} \label{t22}
Suppose $G\alpha$ and $G\beta$ {\em cross}
(i.e., $N(G\alpha,G\beta) \not = 0$).
Just as $\delta$ breaks up into
intersection components of $\alpha$ with $\beta$,
let $\delta' = \beta \cap G\alpha$ give
the analogous disjoint components
$$\matrix{\hfill \delta' =
\gamma_1' \coprod \ldots \coprod \gamma_k' \hfill \cr}$$
with $\gamma_i' = \beta \cap Qx_i^{-1}\alpha$.
Then for all $i$,
$$p|\alpha-\gamma_i|_{aut} + q|\beta-\gamma_i'|_{aut} \leq
p|\alpha|_{aut} + q|\beta|_{aut}.$$
\end{lemma}

\begin{proof} To simplify the notation in the proof below,
we write (just for this proof)
$\|\cdot\|$ for $\|\cdot\|_{aut}$,
$|\cdot|$ for $|\cdot|_{aut}$,
reductive for aut-reductive, etc.

Let $A = \alpha - \gamma_i$ and $B = \beta - \gamma_i'$.
We must show that
$$p|\gamma_i| + q|\gamma_i'| \geq 2\gamma_i.A + 2 \gamma_i'.B.$$
Note that $G\gamma_i=G\gamma_i'$.
Choose coset representatives $y_1, \ldots, y_p$ for
$P$ in $G$ and $z_1, \ldots, z_q$ for $Q$ in $G$.
Then
$$\matrix{
\hfill p|\gamma_i| + q|\gamma_i'| =
\sum_{n=1}^p |\gamma_i| + \sum_{m=1}^q |\gamma_i'| \geq
= 2 |G\gamma_i| = 2 |G\gamma_i'|. \hfill \cr}$$
and
$$2G\gamma_i.A + 2G\gamma_i'.B \leq 2G\gamma_i.(E(\Gamma) - G\gamma_i).$$
So to prove the lemma it suffices to show
$$|G\gamma_i| \geq G\gamma_i.(E(\Gamma) - G\gamma_i),$$
which follows from induction on $|G\gamma_i|$. \end{proof}

Next we review the Pushing and
Shrinking Lemmas of Krstic and Vogtmann
hold in the context of aut-norms and
absolute values.  Unlike the proofs of the
previous two lemma, the proofs for the
next two follow exactly the same lines
as the original proofs by Krstic and Vogtmann
for out-norms and absolute values.
The only way that the new proofs differ
from the old ones is that the new
cardinality conditions for ideal
edges $\alpha_0 \subseteq E_v$ should be verified, namely that:
\begin{itemize}
\item If $v = *$ then $card(\alpha_0) \geq 2$ and
$card(E_v - \alpha_0) \geq 1$.
\item If $v \not = *$ then $card(\alpha_0) \geq 2$ and
$card(E_v - \alpha_0) \geq 2$.
\end{itemize}
As before, it is easily seen from the proofs of
the lemmas that since they hold for both the out-
and aut-norms and absolute values, they also
hold for the tot-norms and absolute values.

\begin{lemma}[Pushing Lemma] \label{t23} 
Let $(\mu,m)$ be a maximally aut-reductive ideal edge
of a reduced pointed marked $G$-graph
with $m \in D(\mu)$.  Let $(\alpha,a)$ be
an aut-reductive ideal edge
containing $m$ which simply crosses $\mu$,
and set $P=stab(\alpha)$.
Then either both $\mu - \alpha$ and $\alpha - \mu$
are aut-reductive or both $\alpha \cup P\mu$ and
$\alpha \cap \mu$ are aut-reductive.
\end{lemma}

\begin{proof} To simplify the notation in the proof below,
we write (just for this proof)
$\|\cdot\|$ for $\|\cdot\|_{aut}$,
$|\cdot|$ for $|\cdot|_{aut}$,
reductive for aut-reductive, etc.

Note that since $m \in \alpha$, $stab(\mu) \leq P$.
As in \cite{[K-V]}, there are four cases depending
upon where $a^{-1}$ and $m^{-1}$ are located.
Since this follows the proof by Krstic and Vogtmann
so closely, the only real detail will be put into the
first case.

\medbreak

\noindent {\em Case 1.} $a^{-1} \not \in G\mu$.
From Lemma \ref{t21},
$$[G:stab(\mu)]|\alpha \cap \mu| + [G:stab(\alpha)]|\alpha \cup P\mu|
\leq [G:stab(\mu)]|\mu| + [G:stab(\alpha)]|\alpha|.$$
Consequently,
$$[G:stab(\mu)](|m|-|\alpha \cap \mu|) +
[G:stab(\alpha)](|a|-|\alpha \cup P\mu|)$$
is greater than or equal to 
$$[G:stab(\mu)](|m|-|\mu|) + [G:stab(\alpha)](|a|-|\alpha|).$$
In other words,
$$red(\alpha \cap \mu,m) + red(\alpha \cup P\mu)
\geq red(\mu,m) + red(\alpha,a).$$

Since $(\mu,m)$ is maximally reductive and $(\alpha,a)$ is
reductive, both of $(\alpha \cup P\mu,a)$
and $(\alpha \cap \mu, m)$ are reductive. As mentioned
above in the discussion preceeding this lemma, we must
verify the cardinality conditions on these two prospective
ideal edges.

First we deal with $(\alpha \cup P\mu,a)$. The
edge $a$ is either bent hyperbolic or elliptic (see
Corollary $4.5$ of \cite{[K-V]}.)  Assume it is bent hyperbolic.
Then as in \cite{[K-V]} we can choose $x \in G$ such that
$xa^{-1} \in E_v-(G\alpha \cup G\mu).$  If $v \not = *$
and $xa^{-1}$ is the only edge in $E_v-(\alpha \cup P\mu)$
then
$$|\alpha \cup P\mu| = |xa^{-1}| = |a^{-1}| = |a|,$$
where the first equality holds because $v \not = *$,
the second is by the $G$-invariance of $|\cdot|$,
and the third follows from our definition of $|\cdot|$
for edges.

In more detail, the first equality
$|\alpha \cup P\mu| = |xa^{-1}|$ holds since
$$\matrix{\hfill E_v = (\alpha \cup P\mu) \coprod \{xa^{-1}\} \hfill \cr}$$
and $v \not = *$.  For a particular coordinate $i$,
both $|\alpha \cup P\mu|_i$ and $|xa^{-1}|_i$ are
measuring the number of times one of the paths $y\alpha_i$
enters $v$ via $\alpha \cup P\mu$ and
leaves via the reverse of $xa^{-1}$ (i.e., $xa$)
\underline{or} enters $v$ via $xa^{-1}$ and leaves it
via the reverse of something in $\alpha \cup P\mu$.
There would be problems if $v=*$ since the
above paths could then enter $v$
and not have to leave it again.

But $|\alpha \cup P\mu| = |a|$ contradicts the fact that
$(\alpha \cup P\mu,a)$ is reductive because
$$[G:stab(\alpha)] (|a| - |\alpha \cup P\mu|) > 0.$$
So if $v \not = *$
then $xa^{-1}$ is not the only edge in $E_v-(\alpha \cup P\mu)$.

For the next possibility, that $a$ is elliptic,
the proof by Krstic and Vogtmann can be used verbatim.

Second we deal with $(\alpha \cap \mu,m)$.  The
set $\alpha \cap \mu$ must contain more than two
edges because it is reductive:
$$[G:stab(\mu)] (|m| - |\alpha \cap \mu|) > 0.$$
The condition on the cardinality of
$E_v - (\alpha \cap \mu)$ is easily
satisfied because $\alpha$ is an ideal
edge and so satisfies the corresponding
condition with $E_v - \alpha$.

\medskip

\noindent {\em Case 2.} $a^{-1} \in G\mu$ and $m^{-1} \in G\alpha$.
Krstic and Vogtmann show that both $(\alpha-\mu,ym^{-1})$
and $(\mu-\alpha,xa^{-1})$ are reductive.

\medskip

\noindent {\em Case 3.} $a^{-1} \in G\mu$,
$m^{-1}  \not \in G\alpha$, and $a \in \mu$.
Both $(\alpha \cup \mu,m)$ and $(\alpha \cap \mu,a)$
are reductive.

\medskip

\noindent {\em Case 4.} $a^{-1} \in G\mu$,
$m^{-1}  \not \in G\alpha$, and $a \not \in \mu$.
Both $(\alpha \cup \mu,m)$ and $(\alpha \cap \mu,m)$
are reductive. \end{proof}

\begin{lemma}[Shrinking Lemma] \label{t24}
 Let $(\mu,m)$ be a maximally aut-reductive \newline
ideal edge
of a reduced pointed marked $G$-graph
with $m \in D(\mu)$.
Let $\alpha$ be
an ideal edge with $N(G\alpha,G\mu) \not = 0$.
Let $\gamma_{i_1}, \ldots, \gamma_{i_k}$ be
the intersection components of $\alpha$ with $\mu$ which
contain no translate of $m$ and let
$\beta = \alpha - \bigcup \gamma_{i_j}$.  Then
$\beta$ or one of the sets $\gamma_{i_j}$ is an
aut-reductive ideal edge.
\end{lemma}

\begin{proof} See verbatim the proof by Krstic and Vogtmann.
If $\alpha_0$ is one of the above
sets, we know it is aut-reductive, and
we want to show it is an ideal
edge, then the cardinality checks
are easy.  The set $\alpha_0$ contains more than
one edge because it is aut-reductive.
Moreover, the cardinality checks on
$E_v - \alpha_0$ follow from similar ones
on $E_v - \alpha$, because
$\alpha_0 \subset \alpha$ for each possibility of
$\alpha_0$. \end{proof}

The following proposition will also be useful in the
next section.

\begin{prop} \label{t25}Let $(\mu,m)$ be a maximally aut-reductive ideal edge
of a reduced pointed marked $G$-graph
with $m \in D(\mu)$.  There is at most one
reductive ideal edge $(\gamma,c)$ at $*$ with
$stab(\gamma)=G$ but where $\gamma$ is not
invertible.  The Whitehead move $(\gamma,c)$ is just
conjugation by $c$, and $\|\gamma\|_{out}=0$.
If $\gamma$ is not compatible with $\mu$, then
$c=m^{-1}$ and $\mu$ is invertible.
\end{prop}

\begin{proof} Since $stab(\gamma)=G$ and $\gamma$ is not invertible,
$E_* - \gamma = \{c^{-1}\}$ must contain just one element.
The Whitehead move
$(\gamma,c)$ consists of first blowing up $\gamma$ and then collapsing
$c$, as illustrated in Figure \ref{fig15}.
The Whitehead move has no effect on the out-norm.  The effect
on the aut-norm can be calculated as follows.
Recall that $F_n = \{\alpha_1, \alpha_2, \ldots \}$.  The
Whitehead move $(\gamma,c)$ conjugates each $\alpha_i$ by
$c$; i.e., each $\alpha_i \mapsto c^{-1}\alpha_ic$ (or
$e(c)^{-1}\alpha_ie(c)$ more accurately, but in the
final graph we can just relabel $e(c)$ as $c$.)

\medskip

\input{fixed.pic}
\begin{figure}[here]
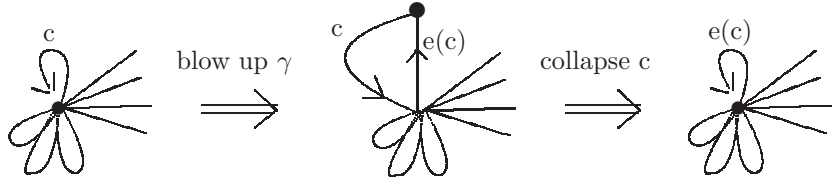

\caption{\label{fig15} The Whitehead move $(\gamma,c)$.}
\end{figure}

Since $\gamma$ is reductive, there exists an $n$ such that
(i) each of $\alpha_1, \ldots, \alpha_{n-1}$ either begins with $c$ 
or ends with $c^{-1}$; and (ii) $\alpha_n$ begins with $c$ and
ends in $c^{-1}$.
Thus for any other $\gamma'=E_* - \{ d^{-1} \}$, $d \not = c$,
$\gamma'$ will increase the length of one of the
$\alpha_1, \ldots, \alpha_n$
and it will not be reductive.  So $(\gamma,c)$ is the only
reductive edge at $*$ with $stab(\gamma)=G$ but where
$\gamma$ is not invertible.

If we further suppose that $\gamma$ is not compatible with $\mu$, then
$c^{-1} \in \mu$ else $\mu \subseteq \gamma$ and they
are compatible.  
As $stab(c^{-1}) = G$, $stab(\mu)=G$.
Since  $\mu$ is reductive and not equal to $\gamma$,
$\mu$ is invertible.  By way of contradiction,
assume $m \in \gamma$.  Then $m \not = c$ else $\gamma$
and $\mu$ are compatible.
We apply the Pushing Lemma to $\gamma$ and $\mu$.
Case $2$ is the relevant case and so
$(\mu-\gamma,a^{-1})$ is reductive.  This
contradicts the fact that $\mu - \gamma = \{a^{-1}\}$
has just one edge in it. So $m \not \in \gamma$
and hence $c=m^{-1}$. \end{proof}
 
\section{The contractibility lemmas} \label{c12}
\arabiceqn

\begin{proof}[Proof of Theorem \ref{tr13}:] 
The space $X_n^G$ deformation retracts to $L_G$.
Following the proof of Theorem $8.1$ by Krstic and Vogtmann
in \cite{[K-V]}, we show that the complex $L_G$ is contractible
by setting
$$L_{<\rho} = \bigcup_{\|\rho'\|_{tot} < \|\rho\|_{tot}} st(\rho')$$
and letting
$$S_\rho= st(\rho) \cap L_{<\rho}.$$
As in \cite{[K-V]}, we show that $S_\rho$ is
contractible when it is non-empty, so that a
transfinite induction argument then yields that
for all $\rho$, all of the components of $L_{<\rho}$
are contractible.  Krstic's work in \cite{[K]} shows that
any two reduced graphs in $L_G$ can be connected by Whitehead moves,
so that $L_G$ is connected.  Thus $L_G$ is contractible
if we can perform the above transfinite induction.

As in \cite{[K-V]}, the first step is to
deformation retract $S_{<\rho}$ to to $S(\mathcal{R})$
by the Poset Lemma (stated in \cite{[K-V]}, deriving from
Quillen in \cite{[Q]}.) We can do this
for the case of $Aut(F_n)$ rather than $Out(F_n)$
without any
significant modifications of the arguments in the previous case.
This is because the Factorization Lemma
and Proposition $6.5$ of $\cite{[K-V]}$ let
us identify $S_{<\rho}$ with the poset of ideal forests
which contain a reductive ideal edge (where we must,
of course, use the newly modified definition of an
ideal forest.)  (The Factorization Lemma gives
a certain isomorphism between forest that does {\em not}
preserve basepoints, but the fact that basepoints are
not preserved is not relevant to Proposition $6.5$.)

After contracting $S_{<\rho}$ to $S(\mathcal{R})$, Krstic
and Vogtmann then use a series of lemmas to
deformation retract from $S(\mathcal{R})$ to $S(\mathcal{C}_1)$,
from there to $S(\mathcal{C}_0)$, and finally
to a point.  We more or less follow this, except there
is an additional
intermediate step where we deformation retract from
$S(\mathcal{C}_1)$ to $S(\mathcal{C}_0^{'})$
and from there to $S(\mathcal{C}_0)$.

The rest of this section will be devoted to proving the
aforementioned series of lemmas which show that
$S(\mathcal{R})$ deformation retracts to a point. \end{proof}

We assume that the maximally reductive ideal edge
$(\mu,m)$ is at the basepoint in all that follows, else the
arguments of Krstic and Vogtmann directly give the
contractibility of $S(\mathcal{R})$.
Moreover, if $(\mu,m) = (\gamma,c)$ where
$\gamma = E_* - \{c^{-1}\}$, then $\mu$ is not
out-reductive at all.  Since $\mu$ is maximally
reductive, Proposition \ref{t25} implies that $\mu$ is the
only reductive ideal edge.  So in this case
$\mathcal{R} = \mathcal{C}_0 = \{\mu\}$ and
$S(\mathcal{R})$ is contractible.  Assume $\mu \not = \gamma$
from now on.

Note the
slight difference in our definition of $\mathcal{C}_1$
from that of \cite{[K-V]},
where here it is phrased to include $\alpha \subset E_v$
which have $stab(\alpha) = stab(v)$, rather than just
invertible $\alpha$.  In other words, from Proposition \ref{t25},
there is at most one reductive ideal edge $(\gamma,c)$
at the basepoint which has $stab(\gamma)=G$ and yet is
not invertible.  This $\gamma$ would be in both $\mathcal{C}_1$
and $\mathcal{C}_0^{'}$.

The next lemma (unlike the ones which follow it) is
essentially the corresponding lemma in \cite{[K-V]} with
minimal modifications.  We repeat their
arguments here for the sake of convenience.

\begin{lemma} \label{t26} The complex $S(\mathcal{R})$ deformation
retracts onto $S(\mathcal{C}_1)$.
\end{lemma}

\begin{proof} Let $\mathcal{C}=\mathcal{C}^\pm$ be a subset of $\mathcal{R}$
which contains $\mathcal{C}_1$.  We show that
$S(\mathcal{C})$ deformation retracts to $S(\mathcal{C}_1)$
by induction on the cardinality of $\mathcal{C} - \mathcal{C}_1$.

Choose $\alpha \in \mathcal{C} - \mathcal{C}_1$ which
satisfies both of:
\begin{enumerate}
\item The cardinality $|\alpha \cap G\mu|$ is minimal
(recall that $\mu$ is the maximally reductive ideal edge.)
\item The ideal edge $\alpha$ is minimal with respect to
property $1.$
\end{enumerate}

Using the Shrinking Lemma 7.4 of \cite{[K-V]}
with $\alpha$ and $\mu$, we obtain
a reductive ideal edge $\alpha_0 \subset \alpha$ which is
compatible with $\mu$.
Let $\gamma_i$ be the intersection components of
$\alpha$ with $\mu$ and index them so that
$m \in \gamma_0$.
Now from the Shrinking Lemma, we
can choose $\alpha_0$ so that it
is either one of the intersection
components $\gamma_{i_j}$ of $\alpha$ with $\mu$ which contain
no translate of $m$, or it is $\alpha - \cup \gamma_{i_j}$.
Because $\alpha \in \mathcal{C}
- \mathcal{C}_1$, $stab(\alpha) \not = G$ and $\alpha$ is neither
invertible nor equal to $\gamma = E_* - \{c^{-1}\}$.

\begin{claim} \label{t27} For every $\beta \in \mathcal{C}$, if $G\beta$
is compatible with $G\alpha$, then $G\beta$ is compatible with
$G\alpha_0$.
\end{claim}

\begin{proof} The three cases are
\begin{enumerate}
\item $G\alpha \subseteq G\beta$.  In this case, $G\alpha_0
\subseteq G\beta$ as $G\alpha_0 \subseteq G\alpha$.
\item $G\alpha \cap G\beta = \void$.
It follows that
$G\alpha_0 \cap G\beta = \void$ since $G\alpha_0 \subseteq G\alpha$.
\item $G\beta \subseteq G\alpha$.
Without loss of generality $\beta \subseteq \alpha$.
If $\beta \not \in \mathcal{C}_1$, then the minimality
conditions on $\alpha$ imply that $\beta=\alpha$, in which case
$\beta$ is clearly compatible with $\alpha_0$.  So
assume $\beta \in \mathcal{C}_1$.
As $G\beta \subseteq G\alpha$, $stab(\beta) \not = G$
since $stab(\alpha) \not = G$.
So either $\beta \in \mathcal{C}_0$ or $m \in G\beta$ and
$N(G\beta,G\mu)=1$.
If $\beta \in \mathcal{C}_0$ then either
$G\beta \subseteq G\mu$ (in which case $\beta$ is in some
$\gamma_i$ and thus compatible with $\alpha_0$),
$G\mu \subseteq G\beta$ (which can not happen as then
$\alpha$ would be in $\mathcal{C}_0$), or
$G\beta \cap G\mu = \void$ (in which case
$\beta \subseteq \alpha - \cup \gamma_{i_j}$ and thus
compatible with $\alpha_0$.)  Finally,
if $m \in G\beta$ and
$N(G\beta,G\mu)=1$ then $\beta \cap G\mu$ is not in
any of the $\gamma_{i_j}$'s as those are the intersection
components of $\alpha$ with $\mu$ that do not
contain a translate of $m$.  In fact, $\beta \cap G\mu$ is
in $\gamma_0$ and
$\beta \subseteq \alpha - \cup \gamma_{i_j}$.  Thus
$\beta$ is compatible with every choice of $\alpha_0$.
\end{enumerate} \end{proof}

Define a poset map
$f : S(\mathcal{C}) \to S(\mathcal{C})$ by sending
an ideal forest $\Phi$ to $\Phi \cup \{G\alpha_0\}$ if
$\Phi$ contains $\alpha$ and
to itself otherwise.  By the Poset Lemma, the
image of $f$ is a deformation retract of $S(\mathcal{C})$,
because $\Phi \subseteq f(\Phi)$
for all $\Phi$.
Define another poset map
$g : f(S(\mathcal{C})) \to f(S(\mathcal{C}))$ by sending
an ideal forest $\Psi$ to $\Psi - \{G\alpha\}$ if
$\Psi$ contains $\alpha$ and
to itself otherwise.  By the Poset Lemma, the
image of $g$ is a deformation retract of $f(S(\mathcal{C}))$,
as $g(\Psi) \subseteq \Psi$ for all $\Psi$.
Hence $S(\mathcal{C})$ deformation
retracts to $S(\mathcal{C}-\{G\alpha\})$, completing
the induction step. \end{proof}

\begin{lemma} \label{t28}   The complex $S(\mathcal{C}_1)$ deformation
retracts onto $S(\mathcal{C}_0^{'})$.
\end{lemma}

\begin{proof} Let $\mathcal{C}$ be a subset of $\mathcal{C}_1$
which contains $\mathcal{C}_0^{'}$.  We show that
$S(\mathcal{C})$ deformation retracts to $S(\mathcal{C}_0^{'})$
by induction on the cardinality of $\mathcal{C} - \mathcal{C}_0^{'}$.

Choose $\alpha \in \mathcal{C} - \mathcal{C}_0^{'}$
such that $m \in \alpha$ and
\begin{enumerate}
\item The cardinality $|\alpha \cap G\mu|$ is minimal.
\item The ideal edge $\alpha$ is minimal with respect to
property $1.$
\end{enumerate}

We apply the Pushing Lemma 7.3 of \cite{[K-V]} 
to get a reductive edge $\alpha_0$
with $\alpha_0=\alpha \cap \mu$ or $\alpha_0 = \alpha - \mu$.
Note that $\alpha_0 \in \mathcal{C}_0$.

\begin{claim} \label{t29} For every $\beta \in \mathcal{C}$, if $G\beta$
is compatible with $G\alpha$, then $G\beta$ is compatible with
$G\alpha_0$.
\end{claim}

\begin{proof} The three cases are
\begin{enumerate}
\item $G\alpha \subseteq G\beta$.  In this case, $G\alpha_0
\subseteq G\beta$ as $G\alpha_0 \subseteq G\alpha$.
\item $G\alpha \cap G\beta = \void$.
It follows that
$G\alpha_0 \cap G\beta = \void$ since $G\alpha_0 \subseteq G\alpha$.
\item $G\beta \subseteq G\alpha$.
Without loss of generality $\beta \subseteq \alpha$.
If $\beta \not \in \mathcal{C}_0^{'}$, then the minimality
conditions on $\alpha$ imply that $\beta=\alpha$, in which case
$\beta$ is clearly compatible with $\alpha_0$.  So
assume $\beta \in \mathcal{C}_0^{'}$.
Since $\alpha \in \mathcal{C} - \mathcal{C}_0^{'}$,
$stab(\alpha) \not = G$.
As $G\beta \subseteq G\alpha$, $stab(\beta) \not = G$
also, which means that $\beta \in \mathcal{C}_0$.
The three ways in which $\beta$ could be compatible with
$\mu$ are:
\begin{itemize}
\item $G\beta \cap G\mu = \void$.
Then $G\beta$ is
disjoint from $G(\alpha \cap \mu)$ and
contained in $G(\alpha - \mu)$.
\item $G\beta \subseteq G\mu$.
Then $G\beta \subseteq G(\alpha \cap \mu)$ and
$G\beta$ is disjoint from $G(\alpha - \mu)$.
\item $G\mu \subseteq G\beta$.
Then $G\mu \subseteq G\alpha$ and so $G\mu$ and $G\alpha$
are compatible, a contradiction.
\end{itemize}
\end{enumerate}\end{proof}

Define a poset map
$f : S(\mathcal{C}) \to S(\mathcal{C})$ by sending
an ideal forest $\Phi$ to $\Phi \cup \{G\alpha_0\}$ if
$\Phi$ contains $\alpha$ and
to itself otherwise.  By the Poset Lemma, the
image of $f$ is a deformation retract of $S(\mathcal{C})$,
because $\Phi \subseteq f(\Phi)$
for all $\Phi$.
Define another poset map
$g : f(S(\mathcal{C})) \to f(S(\mathcal{C}))$ by sending
an ideal forest $\Psi$ to $\Psi - \{G\alpha\}$ if
$\Psi$ contains $\alpha$ and
to itself otherwise.  By the Poset Lemma, the
image of $g$ is a deformation retract of $f(S(\mathcal{C}))$,
as $g(\Psi) \subseteq \Psi$ for all $\Psi$.
Hence $S(\mathcal{C})$ deformation
retracts to $S(\mathcal{C}-\{G\alpha\})$, completing
the induction step. \end{proof}

Now we are left with the task of showing that $S(\mathcal{C}_0^{'})$
deformation retracts to $S(\mathcal{C}_0)$
and from there to a point. The methods used will be analogous
to those in Lemmas \ref{t26} and \ref{t28}, and we will
omit unecessary detail from the remaining proofs.  From
Proposition \ref{t25}, we see that this can be handled in
three separate cases:
\begin{itemize}
\item The ideal edge $\mu$ is invertible and the
reductive ideal edge $\gamma = E_* - \{c^{-1}\}$ is not
compatible with $\mu$.  In this case, the proposition
gives us that $c=m^{-1}$.
\item The ideal edge $\mu$ is invertible and the
reductive ideal edge $\gamma = E_* - \{c^{-1}\}$ is
compatible with $\mu$.
\item The ideal edge $\mu$ is not invertible.
\end{itemize}

\begin{lemma} \label{t30} Suppose $\mu$ is invertible
and $\gamma= E_* - \{m\}$ is reductive.  Then
$S(\mathcal{C}_0^{'})$ is contractible.
\end{lemma}

\begin{proof} We first contract $S(\mathcal{C}_0^{'})$ to
$S(\mathcal{C}_0 \cup \{\gamma\})$.  Let
$\mathcal{C}$ be a subset of $\mathcal{C}_0^{'}$
which contains $\mathcal{C}_0 \cup \{\gamma\}$.
Also assume that if $\alpha \in \mathcal{C}$ is not
pre-compatible with $\mu$, then $\alpha^{-1} \in \mathcal{C}$
also.  We will use induction
on $|\mathcal{C} - (\mathcal{C}_0 \cup \{\gamma\})|$
to show that $S(\mathcal{C})$ deformation retracts to 
$S(\mathcal{C}_0 \cup \{\gamma\})$.

Choose $\alpha \in \mathcal{C} - (\mathcal{C}_0 \cup \{\gamma\})$
such that $m \in \alpha$ and
\begin{enumerate}
\item The cardinality $|\alpha \cap \mu|$ is maximal.
\item The edge $\alpha$ is maximal with respect to property $1.$
\end{enumerate}

There are two main cases, and two subcases in the second case.

\medskip

\noindent {\em Case 1.} $\alpha^{-1}$ is compatible with $\mu$.

Then $\mu \not \subseteq \alpha^{-1}$ because
$m \in \mu$ and $m \in \alpha$. Also,
$\mu \cap \alpha^{-1} \not = \void$ else $\mu \subseteq \alpha$
and $\alpha$ is compatible with $\mu$.  So $\alpha^{-1} \subseteq \mu$.
Let $\alpha_0 = \alpha^{-1}$ and note that $\alpha_0 \in \mathcal{C}$.

For every $\beta$ in $\mathcal{C}$,
if $G\beta$ is compatible with $\alpha$, then $G\beta$
is compatible with $\alpha_0$.  Hence we can replace occurences of
$G\alpha$ in ideal forests with $G\alpha_0$, and retract
$S(\mathcal{C})$ to $S(\mathcal{C} - \{\alpha\})$.

\medskip

\noindent {\em Case 2.} $\alpha^{-1}$ is not compatible with $\mu$.

Since $\alpha$ and $\mu$ cross simply (this is automatic because
$\alpha$ is invertible), the Pushing Lemma applies.
Thus one of the sets $\alpha_0 = \mu - \alpha$ or
$\alpha_0 = \alpha \cup \mu$ is a reductive ideal edge.
As $\mu - \alpha \subseteq \mu$ and $\mu \subseteq \alpha \cup \mu$,
$\alpha_0 \in \mathcal{C}_0$ in either case.

\medskip

\noindent {\em Subcase 1.} $\alpha_0 = \mu - \alpha$.
For every $\beta \in \mathcal{C}$,
if $G\beta$ is pre-compatible with $\alpha$, then
$G\beta$ is compatible with $\alpha_0$. Now replace occurences
of $\alpha$ or $\alpha^{-1}$ with $G\alpha_0$ to retract
$S(\mathcal{C})$ to $S(\mathcal{C} - \{\alpha,\alpha^{-1}\})$.

\medskip

\noindent{Subcase 2.} $\alpha_0 = \alpha \cup \mu$.
Since $m \in \alpha_0$, $\alpha_0 \not = \gamma$.
Accordingly, $\alpha_0$ is invertible and
both $\alpha_0$ and $\alpha_0^{-1}$ are in
$\mathcal{C}_0$.
 
For every $\beta \in \mathcal{C}$,
if $G\beta$ is compatible with $\alpha^{-1}$ then
$G\beta$ is compatible with $\alpha_0^{-1}$.
Substitute $\alpha_0^{-1}$ for $\alpha^{-1}$ to retract
$S(\mathcal{C})$ to
$S(\mathcal{C} - \{\alpha^{-1}\})$.

For every $\beta \in \mathcal{C} - \{\alpha^{-1}\}$,
if $G\beta$ is compatible with $\alpha$ then $G\beta$ is
compatible with $\alpha_0$.
Now substitute $\alpha_0$ for $\alpha$ to retract
$S(\mathcal{C} - \{\alpha^{-1}\})$ to
$S(\mathcal{C} - \{\alpha,\alpha^{-1}\})$.

\medskip

This concludes our argument that
$S(\mathcal{C}_0^{'})$ contracts to
$S(\mathcal{C}_0 \cup \{\gamma\})$.  To eliminate
$\gamma$, note that $\gamma$ is compatible with
$\mu^{-1} \in \mathcal{C}_0$ and verify that
for every $\beta \in \mathcal{C}_0$,
if $G\beta$ is compatible with $\gamma$ then $G\beta$
is compatible with $\mu^{-1}$.
Now replace $\gamma$ with $\mu^{-1}$ to deformation
retract $S(\mathcal{C}_0 \cup \{\gamma\})$
to $S(\mathcal{C})$.

The final step of contracting $S(\mathcal{C}_0)$
to a point is done by adding $\mu$ to all ideal forest
and then removing everything else. \end{proof}

\begin{lemma} \label{t36}  Suppose $\mu$ is invertible
and the reductive $\gamma= E_* - \{c^{-1}\}$ is
compatible with $\mu$.  Then
$S(\mathcal{C}_0^{'})$ is contractible.
\end{lemma}

\begin{proof} The proof of the more complicated case in
Lemma \ref{t30} carries over to this one, with the
exception that the penultimate step of
deformation retracting from
$S(\mathcal{C}_0 \cup \{\gamma\})$
to $S(\mathcal{C})$ is unnecessary, because
$\gamma$ is already compatible with $\mu$.
In addition, various other minor changes need to
be made because $\gamma$ is now compatible with 
$\mu$. \end{proof}

\begin{lemma} \label{t37} Suppose $\mu$ is not invertible.
Then
$S(\mathcal{C}_0^{'})$ is contractible.
\end{lemma}

\begin{proof} As before, let $\gamma = E_* - \{c^{-1}\}$ be the
reductive edge that Proposition \ref{t25} gives us (if it exists).
We know that $\gamma$ is compatible with $\mu$ because
$stab(\mu) \not = G = stab(c^{-1})$.

We first contract $S(\mathcal{C}_0^{'})$ to
$S(\mathcal{C}_0)$.  Let
$\mathcal{C}$ be a subset of $\mathcal{C}_0^{'}$
which contains $\mathcal{C}_0$.
Also assume that if $\alpha \in \mathcal{C}$
and $\alpha$ is invertible, then $\alpha^{-1} \in \mathcal{C}$
also.  We will use induction
on $|\mathcal{C} - \mathcal{C}_0|$
to show that $S(\mathcal{C})$ deformation retracts to 
$S(\mathcal{C}_0)$.

Choose $\alpha \in \mathcal{C} - (\mathcal{C}_0)$
such that $m \in \alpha$ and
\begin{enumerate}
\item The cardinality $|\alpha \cap G\mu|$ is maximal.
\item The edge $\alpha$ is maximal with respect to property $1.$
\end{enumerate}

Since $\alpha$ and $\mu$ cross simply (this is automatic because
$\alpha$ is invertible), the Pushing Lemma applies.
Say $\alpha = (\alpha,a)$.  Neither $a$ nor $a^{-1}$
is in $\mu$ since $stab(a)=G$ and $\mu$ is not invertible
(and not equal to $\gamma$.)
So case 1 of the Pushing Lemma shows that
$\alpha_0 = \alpha \cup G\mu$ is a reductive ideal edge.
As $\mu \subseteq \alpha \cup G\mu$,
$\alpha_0 \in \mathcal{C}_0$.
The ideal edge $\alpha_0$ is not equal to $\gamma$
because it is out-reductive by the proof of the Pushing
Lemma (as both $\alpha$ and $\mu$ are out-reductive.)
So $\alpha_0$ is invertible and
$\alpha_0^{-1} \in \mathcal{C}_0$ also.

As in
Subcase 2, Case 2 of Lemma \ref{t30} above, 
first retract to $S(\mathcal{C} - \{\alpha^{-1}\})$,
then to $S(\mathcal{C} - \{\alpha, \alpha^{-1}\})$, and
finally to $S(\mathcal{C}_0)$. \end{proof}

The sequence of lemmas above concludes our proof of
Theorem \ref{tr13}.

\end{document}